\documentclass[default]{sn-jnl}
\usepackage{graphicx}
\usepackage{multirow}
\usepackage{amsmath,amssymb,amsfonts}
\usepackage{amsthm}
\usepackage{mathrsfs}
\usepackage[title]{appendix}
\usepackage{xcolor}
\usepackage{textcomp}
\usepackage{manyfoot}
\usepackage{booktabs}
\usepackage{algorithm}
\usepackage{algorithmicx}
\usepackage{algpseudocode}
\usepackage{listings}
\theoremstyle{thmstyleone}
\newtheorem{theorem}{Theorem}
\newtheorem{proposition}[theorem]{Proposition}
\newtheorem{lemma}[theorem]{Lemma}
\newtheorem{corollary}[theorem]{Corollary}
\theoremstyle{thmstyletwo}
\newtheorem{example}{Example}
\newtheorem{remark}[theorem]{Remark}
\theoremstyle{thmstylethree}
\newtheorem{definition}[theorem]{Definition}
\begin{document}

\title[Clairaut conformal submersions]{Clairaut Conformal Submersions}
\author[1]{\fnm{Kiran} \sur{Meena}}\email{kirankapishmeena@gmail.com}

\author*[2]{\fnm{Tomasz}\sur{Zawadzki}}\email{tomasz.zawadzki@wmii.uni.lodz.pl}

\affil[1]{\orgdiv{Department of Mathematics}, \orgname{Harish-Chandra Research Institute}, \orgaddress{\street{Chhatnag Road, Jhunsi}, \city{Prayagraj}, \postcode{211019}, \country{India}}}

\affil*[2]{\orgdiv{Faculty of Mathematics and Computer Science}, \orgname{University of Lodz}, \orgaddress{\street{Banacha 22}, \city{Lodz}, \postcode{90-238}, \country{Poland}}}

\abstract{The aim of this paper is to introduce Clairaut conformal submersions between Riemannian manifolds. First, we find necessary and sufficient conditions for conformal submersions to be Clairaut conformal submersions. In particular, we obtain Clairaut relation for geodesics on the total manifolds of conformal submersions, and prove that Clairaut conformal submersions have constant dilation along their fibers, which are totally umbilical, with mean curvature being gradient of a function. Further, we calculate the scalar and Ricci curvatures of the vertical distributions of the total manifolds. Moreover, we find a necessary and sufficient condition for Clairaut conformal submersions to be harmonic. For a Clairaut conformal submersion we find conformal changes of the metric on its domain or image, that give a Clairaut Riemannian submersion, a Clairaut conformal submersion with totally geodesic fibers, or a harmonic Clairaut submersion. Finally, we give two non-trivial examples of Clairaut conformal submersions to illustrate the theory and present a local model of every Clairaut conformal submersion with integrable horizontal distribution.}



\maketitle
\section{Introduction}\label{sec1}
Immersions and submersions are important topics in differential geometry and play a fundamental role in Riemannian geometry. A smooth map $F: (M, g) \to (B, g')$ between Riemannian manifolds is said to be an isometric immersion if the derivative map $F_\ast$ of $F$ is injective and $F$ is an isometry, i.e., $g'(F_\ast X, F_\ast Y) = g(X, Y)~ \text{for}~ X, Y\in \Gamma(TM)$.
On the other hand, Riemannian submersions between Riemannian manifolds were studied by O'Neill \cite{Neill_1966}, Gray \cite{Gray_1967} and Falcitelli \cite{Falcitelli_2008}. A smooth map $F: (M, g) \to (B, g')$ between Riemannian manifolds is said to be a Riemannian submersion if the derivative map $F_\ast$ of $F$ is surjective and satisfies $g'(F_\ast X, F_\ast Y) = g(X, Y)~ \text{for}~ X, Y\in \Gamma(ker F_\ast)^\bot$.
Conformal submersions are a natural generalization of Riemannian submersions \cite{Neill_1966}, as when restricted to vectors orthogonal to their fibers (such vectors are called horizontal) they are conformal maps.
More precisely, a submersion $F:(M,g) \to (B,g')$ 
is said to be a horizontally conformal submersion (or a conformal submersion) if there is a smooth function $\lambda: M \to \mathbb{R}^{+}$ (called dilation) such that
\begin{equation}\label{eqn1.1}
	g'(F_\ast X, F_\ast Y) = \lambda^2 g(X, Y)~ \text{for}~ X, Y \in \Gamma(kerF_\ast)^\bot.
\end{equation}
Obviously, every conformal submersion with $\lambda =1$ is a Riemannian submersion. The function $\lambda^2$ is called the square dilation and it is necessarily positive. In \cite{Gudmundsson_1992}, Gudmundsson obtained the fundamental equations for conformal submersions. The second author studied conformal submersions with totally umbilical, geodesic and minimal fibers in \cite{Zawadzki_2014} and \cite{Zawadzki_2020}. Further, horizontally conformal maps were defined by Fuglede \cite{Fuglede_1978} and Ishihara \cite{Ishihara_1979}, they are useful for characterization of harmonic morphisms \cite{Baird_2003} and have applications in medical imaging (brain imaging) and computer graphics.
A conformal submersion $F:(M,g) \to (B,g')$ is said to be homothetic if the gradient of its dilation $\lambda$ is vertical, i.e., $\mathcal{H} (grad \lambda)= 0$, where $\mathcal{H}$ is the orthogonal projection on the horizontal distribution $(kerF_\ast)^\bot$.

In elementary differential geometry, if $\omega$ is the angle between the velocity vector of a geodesic and a meridian, and $r$ is the distance to the axis of a surface of revolution, then Clairaut's relation (\cite{doCarmo}, \cite{Pressley_2010}) states that $r\sin\omega$ is constant. In \cite{Bishop_1972}, Bishop introduced Clairaut submersions and gave a necessary and sufficient condition for a Riemannian submersion to be Clairaut. Further, Clairaut submersions were studied in \cite{Allison_1996} and \cite{Aso_1991}. Cabrerizo and Fern\'andez obtained the Clairaut relation for geodesics of the Hopf tube in \cite{Cabrerizo_2006}.

The concept of harmonic maps and morphisms is a very useful tool for global analysis and differential geometry. The theory of harmonic maps, developed in \cite{Eells_1964}, is still an active field in differential geometry and has applications to many different areas of mathematics and physics. A map between Riemannian manifolds is called harmonic if the divergence of its differential map vanishes. Harmonic maps between Riemannian manifolds satisfy a system of quasi-linear partial differential equations, therefore to obtain them one would solve partial differential equations on certain manifolds.

In this paper, we study Clairaut conformal submersions between Riemannian manifolds. The paper is organized as follows: in Section \ref{sec2}, we give some basic information about conformal submersions, which is needed in the further parts of this paper. In Section \ref{sec3}, we define Clairaut conformal submersions and obtain necessary and sufficient conditions for a conformal submersion to be Clairaut. For a given Clairaut conformal submersion, we find conformal changes of the metric on its domain or image, that give a Clairaut Riemannian submersion and a Clairaut conformal submersion with totally geodesic fibers. Moreover, we calculate the scalar and Ricci curvatures of the vertical distribution of a Clairaut submersion. In Section \ref{sec4}, we find a necessary and sufficient condition for a Clairaut conformal submersion to be harmonic, and conformal changes of the metric on its domain or image to make a Clairaut conformal submersion harmonic. Finally, in Section \ref{sec5}, we give two non-trivial examples of Clairaut conformal submersions and prove that the domain of every Clairaut conformal submersion with integrable horizontal distribution is locally a warped product.
\section{Preliminaries}\label{sec2}
In this section, we recall the notion of conformal submersion between Riemannian manifolds and give a brief review of basic facts.

Let $F:(M^m,g) \to (B^n,g')$ be a conformal submersion between Riemannian manifolds $(M,g)$ and $(B,g')$, of dimensions $\dim M = m$ and $\dim B =n$. The fiber of $F$ over $q \in B$ is defined as the set $F^{-1}(\{q\})$. The vectors tangent to the fibers form the smooth vertical distribution denoted, at $p \in M$, by $\nu_p$; its orthogonal complement with respect to $g$ is called the horizontal distribution and denoted by $\mathcal{H}_p$. Projections onto the horizontal and vertical distributions are denoted by $\mathcal{H}$ and $\nu$, respectively. A vector field $E$ on $M$ is said to be projectable if there exists a vector field $\tilde{E}$ on $B$ such that $F_{\ast p}(E) = \tilde{E}_{F(p)}$ for all $p \in M$. Then $E$ and $\tilde{E}$ are called $F$-related. For any vector field $\tilde{E}$ on $B$ there exists a unique horizontal vector field $E$ on $M$ such that $E$ and $\tilde{E}$ are $F$-related, then the vector field $E$ is called the horizontal lift of $\tilde{E}$. Horizontal lifts of vector fields from $B$ are also called basic fields on $M$.

The O'Neill tensors $A$ and $T$ were defined in \cite{Neill_1966} as
\begin{equation}\label{eqn2.1}
	A_{X_1} X_2= \mathcal{H} \nabla_{\mathcal{H}X_1} \nu X_2 + \nu \nabla_{\mathcal{H}X_1} \mathcal{H}X_2,
\end{equation}
\begin{equation}\label{eqn2.2}
	T_{X_1} X_2= \mathcal{H} \nabla_{\nu X_1} \nu X_2 + \nu \nabla_{\nu X_1} \mathcal{H}X_2,
\end{equation}
for all $X_1, X_2 \in \Gamma(TM)$, where $\nabla $ is the Levi-Civita connection of $g$. For any $ X_1 \in \Gamma(TM) $, tensors $T_{X_1}$ and $A_{X_1}$ are skew-symmetric operators on $(\Gamma(TM),g)$ reversing the horizontal and the vertical distributions. It is also easy to see that $ T$  is vertical, i.e., $T_{X_1} = T_{\nu X_1} $ and $A$ is horizontal, i.e., $A_{X_1} = A_{\mathcal{H}X_1} $. We note that the tensor field $T$ satisfies $T_U W = T_W U,~ \forall U,W \in \Gamma(kerF_\ast)$. Now, from (\ref{eqn2.1}) and (\ref{eqn2.2}), we have for all $X,Y \in \Gamma(kerF_\ast)^\bot$ and $ U, V \in \Gamma(kerF_\ast)$:
\begin{equation}\label{eqn2.3}
	\nabla_U V = T_U V +  \nu \nabla_U V,
\end{equation}
\begin{equation}\label{eqn2.4}
	\nabla_X U = A_X U +  \nu \nabla_X U,
\end{equation}
\begin{equation}\label{eqn2.5}
	\nabla_X Y = A_X Y +  \mathcal{H} \nabla_X Y.
\end{equation}
A conformal submersion $F$ has totally umbilical fibers  if \cite{Zawadzki_2014}
\begin{equation}\label{eqn2.6}
	T_U V = g(U,V)H ~\text{or}~ T_U X= -g(H, X)U,
\end{equation}
for all $U,V \in \Gamma(kerF_\ast)$ and $X \in \Gamma(kerF_\ast)^\bot$, where $H$ is the mean curvature vector field of the fibers of $F$, given by
\begin{equation}\label{eqn2.7} 
	(m-n) H = \sum\limits_{i=n+1}^{m} T_{U_i} U_i,
\end{equation} 
where $\{U_i\}_{n+1 \leq i \leq m}$ is an orthonormal basis of the fibers of $F$. The horizontal vector field $H$ vanishes if and only if every fiber of $F$ is minimal.
\begin{proposition}\label{prop2.1}
	\cite{Gudmundsson_1992}
	Let $F: (M^m, g) \to (B^n, g')$ be a conformal submersion with dilation $\lambda$. Then
	\begin{equation}\label{eqn2.8}
		A_X Y = \frac{1}{2} \left\{\nu [X, Y] - \lambda^2 g(X, Y) (\nabla_\nu \frac{1}{\lambda^2})\right\}, ~\forall X, Y \in \Gamma(kerF_\ast)^\bot.
	\end{equation}
	Moreover, by (\ref{eqn2.8}), the horizontal space $(kerF_\ast)^\bot$ is totally geodesic if and only if $\lambda$ is constant on $kerF_\ast$.
\end{proposition}

The differential $F_\ast$ of $F$ can be viewed as a section of bundle $Hom(TM,F^{-1}TB)$ $\to M$, where $F^{-1}TB$ is the pullback bundle whose fiber at $p\in M$ is $(F^{-1}TB)_p = T_{F(p)}B$, $p\in M$. The bundle $Hom(TM,F^{-1}TB)$ has a connection $\nabla$ induced from the Levi-Civita connection ${\nabla}^M$ and the pullback connection  ${\nabla}^F$. Then the second fundamental form of $F$ is given by \cite{Nore_1966}
\begin{equation}\label{eqn2.9}
	(\nabla F_\ast) (X,Y) = {\nabla}_X^F F_\ast Y - F_\ast({\nabla}_X^M Y),~\forall X,Y \in \Gamma(TM),
\end{equation} 
or
\begin{equation}\label{eqn2.10}
	(\nabla F_\ast) (X,Y) = {\nabla}_{F_\ast X}^B F_\ast Y - F_\ast({\nabla}_X^M Y),~\forall X,Y \in \Gamma(TM),
\end{equation}
where $\nabla^B$ is the Levi-Civita connection on $B$. Note that for the sake of simplicity we can write $\nabla^M$ as $\nabla$.
\begin{lemma}\label{lem2.1}
	\cite{Gudmundsson_1992}
	Let $F: (M^m, g) \to (B^n, g')$ be a conformal submersion. Then
	\begin{equation*}
		\begin{array}{ll}
			F_\ast(\mathcal{H} \nabla_X Y) = \nabla_{F_\ast X}^B F_\ast Y + \frac{\lambda^2}{2} \left\{X (\frac{1}{\lambda^2})F_\ast Y + Y (\frac{1}{\lambda^2})F_\ast X - g(X, Y) F_\ast (grad_\mathcal{H}(\frac{1}{\lambda^2}))\right\},
		\end{array}
	\end{equation*}
	for all basic vector fields $X, Y$, where $\nabla$ is the Levi-Civita connection on $M$.
\end{lemma}
\noindent Now, using (\ref{eqn2.9}) in Lemma \ref{lem2.1}, we get
\begin{equation}\label{eqn2.11}
	\begin{array}{ll}
		(\nabla F_\ast) (X,Y) =-\frac{\lambda^2}{2} \left\{X (\frac{1}{\lambda^2})\tilde{Y} + Y (\frac{1}{\lambda^2})\tilde{X} - g(X, Y) F_\ast (grad_\mathcal{H}\frac{1}{\lambda^2}))\right\},
	\end{array}
\end{equation}
where $X$ and $Y$ are horizontal lifts of $\tilde{X}$ and $\tilde{Y}$, respectively.

Now we recall the definitions of gradient, divergence and Laplacian \cite{Sahin_2017}. Let $f\in \mathcal{F}(M)$, then the gradient of $f$, denoted by $\nabla f$ or $gradf$, is given by
\begin{equation}\label{eqn2.12} 
	g (gradf, X) = X(f),~\forall X \in \Gamma(TM).
\end{equation}
The divergence of $X$, denoted by $div(X)$, is given by
\begin{equation}\label{eqn2.13}
	div(X)  = \sum_{k=1}^{m} g(\nabla_{e_k} X, e_k), ~\forall X\in \Gamma(TM),
\end{equation}
where $\{e_k\}_{1\leq k \leq m}$ is an orthonormal basis of $T_pM$.
The Laplacian of $f$, denoted by $\Delta f$, is given by
\begin{equation}\label{eqn2.14}
	\Delta f = div(\nabla f).	
\end{equation}
\section{Clairaut Conformal Submersions}\label{sec3}
In this section, we find some interesting results which are useful to investigate the geometry of Clairaut conformal submersions.

The Clairaut condition for Riemannian submersions was originally defined by Bishop in \cite{Bishop_1972}. Further, Clairaut condition was defined for  Riemannian maps in \cite{Sahin_2017b} and \cite{Meena}. Now, we are going to define the Clairaut condition for conformal submersions. 

\begin{definition}\label{def3.1}
	A conformal submersion $F: (M^m, g)\to (B^n, g')$ between Riemannian manifolds with dilation $\lambda$ is said to be a Clairaut conformal submersion if there is a function $r: M \to \mathbb{R}^{+}$ such that for every geodesic $\alpha$ on $M$, the function 
	\linebreak 
	$(r \circ \alpha) \sin \omega(t)$ is constant along $\alpha$, where for all $t$, $\omega(t)$ is the angle between $\dot{\alpha}(t)$ and the horizontal space at $\alpha(t)$.
\end{definition}

We note that the above definition does not depend on the parametrization of $\alpha$. Indeed, $\omega(t) \in [0,\frac{\pi}{2}]$ is the angle between $\dot{\alpha}(t)$ and $\mathcal{H} \dot{\alpha}(t)$, so it does not depend on the length of $\dot{\alpha}(t)$ or the direction in which $\alpha$ is travelled, only on a point on the geodesic. Hence, from now on we will only consider geodesics with constant speed $\|\dot{\alpha}(t)\|$.

We know that with help of Clairaut's theorem we can find all geodesics on a surface of revolution, and the notion of Clairaut conformal submersion comes from geodesic curves on a surface of revolution. Therefore, motivated by \cite{Sahin_2017b}, we are going to find necessary and sufficient conditions for a curve on the total space, and its projection onto the base space, to be geodesics.
\begin{proposition}\label{prop3.1}
	Let $F: (M^m, g)\to (B^n, g')$ be a conformal submersion. Let $\alpha: I \to M$ be a regular curve on $M$, parametrized by arc length, such that $U(t) = \nu \dot{\alpha}(t)$ and $X(t)= \mathcal{H} \dot{\alpha}(t)$. Then $\alpha$ is geodesic on $M$ if and only if 
	\begin{equation}\label{eqn3.1}
		A_X X + \nu \nabla_X U + T_U X + {\nu \nabla_U U}=0,
	\end{equation}
	and
	\begin{equation}\label{eqn3.2}
		\mathcal{H}\nabla_X X + 2 A_X U + T_U U =0.
	\end{equation}
\end{proposition}
\begin{proof}
	Let $\alpha: I \to M$ be a regular curve on $M$ and $\dot{\alpha} = X(t) + U(t)$, where  $U(t) \in \Gamma(kerF_\ast)$ and $X(t) \in \Gamma(kerF_\ast)^\bot$ is the basic vector field along $\alpha$. Then 
	\begin{equation*}
		\nabla_{\dot{\alpha}} \dot{\alpha} = \nabla_{X + U} X+U= \nabla_X X + \nabla_X U +\nabla_U X +\nabla_U U,
	\end{equation*}
	which implies
	\begin{equation*}
		\nabla_{\dot{\alpha}} \dot{\alpha} = \nu \nabla_X X + \nu \nabla_X U + \nu \nabla_U X + \nu \nabla_U U + \mathcal{H} \nabla_X X + \mathcal{H} \nabla_X U +\mathcal{H} \nabla_U X +\mathcal{H} \nabla_U U.
	\end{equation*}
	Using (\ref{eqn2.2}), (\ref{eqn2.3}), (\ref{eqn2.4}) and (\ref{eqn2.5}) in above equation, we get
	\begin{equation}\label{eqn3.3}
		\nabla_{\dot{\alpha}} \dot{\alpha} = A_X X + \nu \nabla_X U + T_U X + \nu \nabla_U U + \mathcal{H} \nabla_X X +2 A_X U +T_U U.
	\end{equation}
	We know that $\alpha$ is geodesic curve if and only if $\nabla_{\dot{\alpha}} \dot{\alpha} =0$. Then proof follows by (\ref{eqn3.3}).
\end{proof}

\begin{theorem}\label{thm3.1}
	Let $F: (M^m, g)\to (B^n, g')$ be a conformal submersion with connected fibers and dilation $\lambda$. Then $F$ is a Clairaut conformal submersion with $r = e^f$ if and only if $\nabla f$ is horizontal, fibers of $F$ are totally umbilical with $H = -  \nabla f$ and $\lambda$ is constant along the fibers of $F$.
\end{theorem}
\begin{proof}
	First we prove that $F$ is a Clairaut conformal submersion with $r = e^f$ if and only if for any geodesic $\alpha: I \rightarrow M$ with $U(t)= \nu \dot{\alpha}(t)$ and $X(t)=  \mathcal{H} \dot{\alpha}(t)$, $t \in I \subset \mathbb{R}$, the equation
	\begin{equation}\label{eqn3.5a}
		\begin{array}{ll}
			g(U(t), U(t)) g(\dot{\alpha}(t), (\nabla f)_{\alpha(t)} ) \\ + \frac{\lambda^2}{2}g\left(\nabla_\nu \frac{1}{\lambda^2}, U(t) \right) g(X(t), X(t)) + g((T_U U)(t), X(t)) = 0,
		\end{array}
	\end{equation}
	is satisfied. To prove this, let $\alpha$ be a geodesic on $M$ with $\dot{\alpha}(t) = X(t) + U(t)$ and let $\omega(t) \in [0,\frac{\pi}{2}]$ denote the angle between $\dot{\alpha}(t)$ and $X(t)$. For all $t \in I$, for which $\alpha(t)$ is horizontal, we have $U(t) = 0$ and (\ref{eqn3.5a}) is satisfied, also $\sin \omega (t) = 0$ and $(r(\alpha(t))) \sin \omega(t)$ identically vanishes for any function $r= e^f$ on $M$. Therefore, the statement holds trivially in this case. 
	Now, we consider a non-horizontal geodesic $\alpha$, then we have $\sin \omega (t) \neq 0$ on an open subset $J$ of $I$. By the previous argument,
	$(r \circ \alpha) \sin \omega (t)$ is constant on $I \setminus J$ if and only if (\ref{eqn3.5a}) holds there, and since
	\linebreak 
	$(r \circ \alpha) \sin \omega (t)$ is continuous with respect to $t$, it is enough to prove that it is constant on $J$ if and only if (\ref{eqn3.5a}) holds on $J$. 
	Since $\alpha$ is a geodesic, its speed is constant $a = \|\dot{\alpha}\|^2$ (say).
	Now
	\begin{equation*}
		\cos^2 \omega(t) = \frac{g(\dot{\alpha}(t), X(t)) g(\dot{\alpha}(t), X(t))}{\|\dot{\alpha}(t)\|^2 \|X(t)\|^2},
	\end{equation*}
	which implies
	\begin{equation}\label{eqn3.6}
		g(X, X) = a \cos^2 \omega(t).
	\end{equation}
	Similarly, we can get
	\begin{equation}\label{eqn3.7}
		g(U, U) = a \sin^2 \omega(t).
	\end{equation}
	Differentiating (\ref{eqn3.7}), we get
	\begin{equation}\label{eqn3.8}
		\frac{d}{dt} g(U, U) = 2g(\nabla_{\dot{\alpha}}U, U) = 2g(\nu \nabla_{U}U + \nu\nabla_{X}U, U).
	\end{equation}
	Using (\ref{eqn3.1}) in (\ref{eqn3.8}), we get
	\begin{equation}\label{eqn3.9}
		\frac{d}{dt} g(U, U) = -2g(A_X X + T_U X, U).
	\end{equation}
	On the other hand, by (\ref{eqn3.7}), we have
	\begin{equation}\label{eqn3.10}
		\frac{d}{dt} g(U, U) = 2 a \sin \omega (t) \cos \omega (t) \frac{d \omega } {dt}.
	\end{equation}
	By (\ref{eqn3.9}) and (\ref{eqn3.10}), we get
	\begin{equation}\label{eqn3.11}
		g(A_X X + T_U X, U) =- a \sin \omega (t) \cos \omega (t) \frac{d \omega } {dt}.
	\end{equation}
	Moreover, $F$ is a Clairaut conformal submersion with $r=e^f \iff \frac{d}{dt} (e^{f\circ\alpha} \sin\omega) = 0 \iff \cos\omega \frac{d\omega}{dt} + \sin\omega \frac{df}{dt}= 0$. Multiplying by non-zero factor $a\sin\omega$, we can write the last equality as
	\begin{equation}\label{eqn3.12}
		- a\cos\omega \sin\omega \frac {d\omega}{dt}= a \sin^2\omega \frac {df}{dt}.
	\end{equation}
	By (\ref{eqn3.7}), (\ref{eqn3.12}) is equivalent to
	\begin{equation}\label{eqn3.13}
		g(U, U)\frac {df}{dt} = - a\cos\omega \sin\omega \frac {d\omega}{dt}.
	\end{equation}
	By (\ref{eqn3.11}), we can write (\ref{eqn3.13}) as
	\begin{equation*}
		g(U, U)\frac {df}{dt} = g(A_X X + T_U X, U),
	\end{equation*}
	which is equivalent to
	\begin{equation}\label{eqn3.14}
		-g(T_U U, X) = g(U, U)\frac {df}{dt} - g(A_X X , U).
	\end{equation}
	Using (\ref{eqn2.8}) in (\ref{eqn3.14}), we get
	\begin{equation*}
		-g(T_U U, X) = g(U, U)\frac {df}{dt} + \frac{\lambda^2}{2} g(X, X) g\left(U, \nabla_\nu \frac{1}{\lambda^2}\right),
	\end{equation*}
	which is equivalent to
	\begin{equation*}
		g(T_U U, X) = -g(U, U) g(\dot{\alpha}, \nabla f) - \frac{\lambda^2}{2} g(X, X) g\left(U, \nabla_\nu \frac{1}{\lambda^2}\right),
	\end{equation*}
	which is same as (\ref{eqn3.5a}). 
	
	If we consider any geodesic $\alpha$ on $M$ with initial vertical tangent vector i.e. $X=0$, then by above equation $g(U, U) g(U, \nabla f) = 0 \implies U(f) = 0$. Thus $f$ is constant on any fiber, as fibers are connected. Therefore $gradf$ turns out to be horizontal. Now, putting $\dot{\alpha} = X+ U$ in above equation, we get
	\begin{equation}\label{eqn3.15}
		g(T_U U, X) =- g(U, U) g(X, \nabla f) - \frac{\lambda^2}{2} g(X, X) g\left(U, \nabla_\nu \frac{1}{\lambda^2}\right).
	\end{equation}
	On the other hand, let us consider the geodesic with initial tangent vector $U + cX$, for arbitrary $c \neq 0$. Then we get from (\ref{eqn3.15}):
	\begin{equation}\label{eqn3.15c}
		c g(T_U U, X) =- c g(U, U) g(X, \nabla f) - c^2 \frac{\lambda^2}{2} g(X, X) g\left(U, \nabla_\nu \frac{1}{\lambda^2}\right).
	\end{equation}
	Taking the difference of (\ref{eqn3.15c}) and (\ref{eqn3.15}) yields
	\[
	(c-1) \left( g(T_U U, X) +  g(U, U) g(X, \nabla f) + (c+1) \frac{\lambda^2}{2} g(X, X) g\left(U, \nabla_\nu \frac{1}{\lambda^2}\right) \right) = 0,
	\]
	which must hold for all $c \neq 0$. Hence, we obtain a pair of equations:
	\begin{equation}\label{fibersumbilical}
		g(T_U U, X) +  g(U, U) g(X, \nabla f) =0
	\end{equation}
	and
	\begin{equation}\label{dilationconstalongfibers}
		\frac{\lambda^2}{2} g(X, X) g\left(U, \nabla_\nu \frac{1}{\lambda^2}\right) =0,
	\end{equation}
	for all vertical $U$ and horizontal $X$. 
	
	From (\ref{fibersumbilical}) it follows that the fibers are totally umbilical with mean curvature vector $-\mathcal{H} \nabla f$. Indeed, for one-dimensional fibers we obtain this immediately from (\ref{fibersumbilical}); if the fibers are of dimension $2$ or greater, let $V,W$ be orthogonal vertical vectors. Then from $g(W,V)=0$, $T_W V = T_V W$ and (\ref{fibersumbilical}) for $U=W+V$ and $U=W-V$, respectively, we obtain:
	\begin{equation}\label{fibersumbilicalUplusV}
		g(T_W W, X) + g(T_V V, X) + 2 g(T_W V , X) + g(W,W) g(X, \nabla f) +  g(V, V) g(X, \nabla f) =0 
	\end{equation}
	and
	\begin{equation}\label{fibersumbilicalUminusV}
		g(T_W W, X) + g(T_V V, X) - 2 g(T_W V , X) + g(W,W) g(X, \nabla f) +  g(V, V) g(X, \nabla f) =0,
	\end{equation}
	and taking the difference of (\ref{fibersumbilicalUplusV}) and (\ref{fibersumbilicalUminusV}) we obtain  $T_W V =0$ for all orthogonal $W$ and $V$. Therefore, we can write for all vertical $U,V$ and horizontal $X$:
	\begin{equation}\label{fibersumbilicalUV}
		g(T_U V, X) = -  g(U, V) g(X, \nabla f),
	\end{equation}
	and hence the fibers of $F$ are totally umbilical with mean curvature $H = - \mathcal{H} \nabla f = -\nabla f$.
	
	On the other hand, from (\ref{dilationconstalongfibers}) we obtain $g\left(U, \nabla_\nu \frac{1}{\lambda^2}\right) =0$ for all vertical $U$, and hence $\nabla_\nu \frac{1}{\lambda^2}  =0$, i.e., $\lambda$ is constant along fibers.
\end{proof}

\begin{remark}
	From Definition \ref{def3.1} it follows that on the domain of a Clairaut conformal submersion a geodesic horizontal at one point is horizontal everywhere, which is a known property of Riemannian submersions \cite{Neill_1966}. However, as evident from its proof, Theorem \ref{thm3.1} remains true if Definition \ref{def3.1} is formulated only for nowhere horizontal geodesics.
\end{remark}

Given a conformal submersion $F: (M, g)\to (B, g')$ with totally umbilical fibers, the mean curvature of which is gradient of some function, we can obtain another mapping with all these properties, by conformally changing metric $g$ (by a function that is constant along the fibers) or $g'$. In what follows, we use these transformations to relate Clairaut conformal submersions to: Clairaut Riemannian submersions, Clairaut conformal submersions with totally geodesic fibers, and, in Section~\ref{sec4}, harmonic Clairaut conformal submersions.

Theorem \ref{thm3.1} states the same conditions for geometry of fibers, $H$ and $f$, as those obtained by considering Riemannian submersion in \cite{Bishop_1972}. Therefore, we can obtain the following:
\begin{corollary}
	Let $F: (M^m, g)\to (B^n, g')$ be a Clairaut conformal submersion with connected fibers, dilation $\lambda$ and $r = e^f$. Let $\psi$ be the function on $B$ such that for all $x \in M$ we have $\psi (F(x)) = \frac{1}{  \lambda(x) }$. Then $F: (M, g)\to (B, \psi^2 g')$ is a Clairaut Riemannian submersion with $r = e^f$.
\end{corollary}
\begin{proof}
	From Theorem \ref{thm3.1} it follows that $\psi$ is well-defined, as for all $x \in M$ if $F(x) = F(y)$ then $\lambda (x) = \lambda (y)$. From (\ref{eqn1.1}) it follows that
	\[
	(\psi(F(p)) )^2 g'(F_\ast X, F_\ast Y) = (\psi(F(p)) )^2 \lambda^2(p) g(X, Y) =  g(X, Y)~ \text{for}~ X, Y \in \mathcal{H}_p\]
	and hence $F: (M, g)\to (B, \psi^2 g')$ is a Riemannian submersion. As the fibers of $F$ are totally umbilical with mean curvature $H = - \mathcal{H} \nabla f = -\nabla f$, from \cite{Bishop_1972} it follows that $F: (M, g)\to (B, \psi^2 g')$ is a Clairaut submersion with $r = e^f$.
\end{proof}

\begin{lemma} \label{lemmaconformalchangeonM}
	Let $F: (M^m, g)\to (B^n, g')$ be a Clairaut conformal submersion with connected fibers, dilation $\lambda$ and $r = e^f$. Let $\phi$ be a positive function on $M$, such that $\nu \nabla \phi=0$. 
	Then $F : (M, \phi^2  g )\to (B, g')$ is a Clairaut conformal submersion with dilation $\frac{\lambda}{\phi}$ and $r = \phi e^f$.
\end{lemma}
\begin{proof}
	Let  ${\tilde F} : (M, \phi^2 g  )\to (B, g')$ be such that ${\tilde F}(x) = F(x)$. Then the fibers of $F$ and ${\tilde F}$ coincide, the horizontal distributions of $F$ and ${\tilde F}$ coincide, and from (\ref{eqn1.1}) it follows that ${\tilde F}$ is a conformal submersion with dilation $\frac{\lambda}{\phi}$. Since totally umbilical submanifolds remain totally umbilical after a conformal change of metric, the fibers of ${\tilde F}$ are totally umbilical on $(M, \phi^2 g  )$. Let $\{ U_i\}_{n+1 \leq i \leq m}$ be a local $g$-orthonormal basis of the fibers of $F$, then $\{ \frac{1}{\phi} {\tilde U}_i\}_{n+1 \leq i \leq m}$ is a local $( \phi^2 g)$-orthonormal basis of the fibers of ${\tilde F}$. Let ${\tilde \nabla}$ be the Levi-Civita connection of $( \phi^2 g)$, then
	\begin{align*}
		\mathcal{H} {\tilde \nabla}_{{\tilde U}_i} {\tilde U}_i &=  \mathcal{H} \nabla_{ \frac{1}{\phi} U_i} ( \frac{1}{\phi} U_i) - \frac{1}{\phi^2}   g(U_i , U_i) \mathcal{H} \nabla \log \phi \\
		&=   \frac{1}{\phi^2} ( - \nabla f - \nabla \log \phi). 
	\end{align*}
	It follows that the mean curvature of the fibers of ${\tilde F}$ on $(M, \phi^2 g  )$ is 
	\[
	{\tilde H} = - \frac{1}{\phi^2} \nabla ( f + \log \phi). 
	\]
	Since for all $X \in \Gamma(TM)$ and all differentiable functions $\psi$ on $M$ we have $g(X, \nabla \psi)= X(\psi) = \phi^2 g( {\tilde \nabla} \psi , X)$, it follows that 
	\begin{equation}\label{eqHconformallychanged}
		{\tilde H} = - {\tilde \nabla} ( f + \log \phi) = - {\tilde \nabla} \log (\phi e^f).
	\end{equation}
\end{proof}

\begin{corollary}\label{cor3.3}
	Let $F: (M^m, g)\to (B^n, g')$ be a Clairaut conformal submersion with connected fibers, dilation $\lambda$ and $r = e^f$. Then $F: (M, \lambda^2 g)\to (B, g')$ is a Clairaut Riemannian submersion with $r = \lambda e^{f }$.
\end{corollary}
\begin{proof}
	From  (\ref{eqn1.1}) it follows that $F: (M, \lambda^2 g)\to (B, g')$ is a Riemannian submersion, and since $\nu \nabla \lambda =0$ by Theorem \ref{thm3.1}, from Lemma \ref{lemmaconformalchangeonM} it follows that $F: (M, \lambda^2 g)\to (B, g')$ is a Clairaut submersion.
\end{proof}

\begin{corollary}
	Let $F: (M^m, g)\to (B^n, g')$ be a Clairaut conformal submersion with connected fibers, dilation $\lambda$ and $r = e^f$.  Then $F: (M, e^{-2f}g)\to (B, g')$ is a Clairaut conformal submersion with totally geodesic fibers and $r=1$.
\end{corollary}
\begin{proof}
	Since $\nu \nabla f =0$ by Theorem \ref{thm3.1}, from Lemma \ref{lemmaconformalchangeonM} it follows that $F: (M, \lambda^2 g)\to (B, g')$ is a Clairaut conformal submersion, from (\ref{eqHconformallychanged}) it follows that the mean curvature of its totally umbilical fibers vanishes, and hence those fibers are totally geodesic.
\end{proof}

\begin{theorem}\label{thm3.3}
	Let $F: (M^m, g) \to (B^n, g')$ be a Clairaut conformal submersion with $r = e^f$ and dilation $\lambda$. Then 
	\begin{align*}
		K_\nu& = \hat{K} - (m-n)(m-n-1) \| \nabla f \|^2,
	\end{align*}
	where $K_\nu$ and $\hat{K}$ denote the scalar curvatures of the vertical distribution and the fibers, respectively.
\end{theorem}
\begin{proof}
	Let $g$ and $\tilde{g}$ be two Riemannian metrics on $M$ which are conformally related to each other. 
	Let $F: (M,g) \to (B, g')$ be a Clairaut conformal submersion with $r=e^f$, then by Theorem \ref{thm3.1}, we have
	\begin{equation}\label{eqn3.30}
		T_U V =- g(U, V) \nabla f; ~\nu (\nabla \lambda) =0,
	\end{equation}
	where $U, V \in \Gamma(kerF_\ast)$. Also, $T$ denotes the O'Neill tensor for the Levi-Civita connection $\nabla$ of $g$ and $\lambda$ is a positive function on $M$. In addition, by Corollary \ref{cor3.3} for $\tilde{g} = \lambda^2 g$ the submersion $\tilde{F}: (M,\tilde{g}) \to (B, g')$ such that $\tilde{F}(p) = F(p)$ for all $p \in M$ is a Clairaut Riemannian submersion with $r = \lambda e^f$. Easily, we see that $F$ and $\tilde{F}$ have the same vertical and horizontal distributions. Then for any vertical vectors $U$ and $V$, we have \cite{Bishop_1972}
	\begin{equation}\label{eqn3.29}
		\tilde{T}_U V = \tilde{g}(U, V) \tilde{H}= - g(U, V) \nabla (f + \log \lambda),
	\end{equation}
	where $\tilde{H} = - {{\tilde \nabla} \log (\lambda e^f)}$ is the mean curvature vector field of the fibers of $\tilde{F}$ and $\tilde{T}$ denotes the O'Neill	tensor for the Levi-Civita connection $\tilde{\nabla}$ of $\tilde{g}$. Now, using (\ref{eqn3.30}), (\ref{eqn3.29}) and (\ref{eqn2.12}), we get
	\begin{equation}\label{eqn3.31}
		\tilde{T}_U V = T_U V - g(U, V) \nabla ( \log \lambda).
	\end{equation}
	Let $R$ and $\tilde{R}$ denote the Riemannian curvature tensors with respect to $\nabla$ and $\tilde{\nabla}$, respectively. Then 
	\begin{align*}
		\tilde{sec}(U, V)= \frac{\tilde{g}(\tilde{R}(U,V)V, U)}{\tilde{g}(U, U) \tilde{g}(V, V) - \tilde{g}(U, V)^2 } = \frac{1}{\lambda^2} g(\tilde{R}(U,V)V, U).
	\end{align*}
	Using (\cite{Gudmundsson_1992}, p. 14) in the above equation, we get
	\begin{align*}
		\lambda^{2} \tilde{sec}(U, V)&= sec(U, V) + \frac{\lambda^2}{2} \left( Hess\frac{1}{\lambda^2} (U, U) + Hess\frac{1}{\lambda^2} (V, V)\right) \nonumber \\& - \frac{\lambda^4}{4} \mid grad \frac{1}{\lambda^2}\mid^2 - \frac{\lambda^4}{4}\left( \left( U\frac{1}{\lambda^2}\right)^2 + \left( V \frac{1}{\lambda^2}\right)^2 \right),
	\end{align*} 
	which implies
	\begin{align}\label{eqn3.32}
		\lambda^{2} \tilde{sec}(U, V)& = sec(U, V) + \frac{\lambda^2}{2} \left( -g(T_U U, \nabla_{\mathcal{H}} \frac{1}{\lambda^2}) - g(T_V V, \nabla_{\mathcal{H}} \frac{1}{\lambda^2})\right)  - \frac{\lambda^4}{4} \mid \nabla_{\mathcal{H}} \frac{1}{\lambda^2}\mid^2
	\end{align} 
	where $sec(U, V)$ and $\tilde{sec}(U, V)$ denote the sectional curvatures of the plane spanned by (orthonormal with respect to $g$) vectors $U$ and $V$, on $(M, g)$ and $(M, \tilde{g})$, respectively. In addition, $Hess\frac{1}{\lambda^2} (U, U)= g(\nabla_U \nabla \frac{1}{\lambda^2},U)$. Since $\tilde{F}: (M,\tilde{g}) \to (B, g')$ is a Clairaut Riemannian submersion, for non-zero orthogonal (with respect to $\tilde{g}$) vertical vectors at $p \in M$, we have (\cite{Falcitelli_2008}, p. 12)
	\begin{equation*}
		\tilde{g} (\tilde{R}(U,V)V, U) = \tilde{g} (\hat{\tilde{R}}(U, V)V, U) + \tilde{g}(\tilde{T}_U V, \tilde{T}_V U) - \tilde{g}(\tilde{T}_V V, \tilde{T}_U U),
	\end{equation*}
	where $\hat{\tilde{R}}$ is the Riemannian curvature tensor of the fibers of $\tilde{F}$. Using (\ref{eqn3.29}) in the above equation, we get
	\begin{equation}\label{eqn3.33}
		\tilde{sec}(U, V) = \hat{\tilde{sec}}(U, V) - \frac{1}{\lambda^2} \mid \nabla (f + \log \lambda )\mid^2,
	\end{equation}
	where $\hat{\tilde{sec}}$ denotes the sectional curvature of the fibers of $\tilde{F}$. We know that 
	\begin{align*}
		\hat{\tilde{sec}}(U, V)= \frac{\tilde{g}(\hat{\tilde{R}}(U,V)V, U)}{\tilde{g}(U, U) \tilde{g}(V, V)- \tilde{g}(U, V)^2 } = \frac{1}{\lambda^2} g(\hat{\tilde{R}}(U,V)V, U).
	\end{align*}
	Using (\cite{Gudmundsson_1992}, p. 14) in the above equation, we get
	\begin{align*}
		\lambda^{2} \hat{\tilde{sec}}(U, V)&= \hat{s}ec(U, V) + \frac{\lambda^2}{2} \left( g(\nabla_U  \nabla_{\nu} \frac{1}{\lambda^2}, U) + g(\nabla_V  \nabla_{\nu} \frac{1}{\lambda^2}, V) \right) \nonumber \\& - \frac{\lambda^4}{4} \mid grad_{\nu} \frac{1}{\lambda^2}\mid^2 - \frac{\lambda^4}{4}\left( \left( U\frac{1}{\lambda^2}\right)^2 + \left( V\frac{1}{\lambda^2}\right)^2 \right),
	\end{align*} 
	which becomes
	\begin{align}\label{eqn3.34}
		\lambda^{2} \hat{\tilde{sec}}(U, V)= \hat{s}ec(U, V).
	\end{align}
	By using (\ref{eqn3.32}), (\ref{eqn3.33}) and (\ref{eqn3.34}), we get
	\begin{align}\label{eqn3.35}
		\frac{1}{\lambda^2}	\hat{s}ec(U, V) - \frac{1}{\lambda^2} \mid \nabla (f + \log \lambda )\mid^2&= \frac{1}{\lambda^2} \Big( sec(U, V) - \frac{\lambda^4}{4} \mid \nabla_{\mathcal{H}} \frac{1}{\lambda^2}\mid^2 \nonumber\\& - \frac{\lambda^2}{2} \left( g(T_U  U, \nabla \frac{1}{\lambda^2}) + g(T_V V, \nabla \frac{1}{\lambda^2}) \right) \Big).
	\end{align}
	Using (\ref{eqn3.30}) in (\ref{eqn3.35}), we get
	\begin{align}\label{eqn3.37}
		sec(U, V)& = \hat{s}ec(U, V) - \mid \nabla (f + \log \lambda )\mid^2 + \frac{\lambda^4}{4} \mid grad_{\mathcal{H}} \frac{1}{\lambda^2}\mid^2 -\lambda^2g\left(\nabla f, \nabla_{\mathcal{H}} \frac{1}{\lambda^2}\right) \nonumber \\
		&= \hat{s}ec(U, V) - \mid \nabla (f + \log \lambda )\mid^2 + \mid \nabla \log \lambda \mid^2 + 2 g\left(\nabla f, \nabla \log \lambda \right) \nonumber \\
		&= \hat{s}ec(U, V) - \mid \nabla f \mid^2 .
	\end{align}
	Note that the same result follows directly from \cite[Theorem 2.2.3(1)]{Gudmundsson_1992}.
	
	Taking trace of (\ref{eqn3.37}), we get
	\begin{align*}
		\sum\limits_{i=1}^{m-n} \sum \limits_{j=1, j \neq i}^{m-n} sec(U_i, U_j) &= \sum\limits_{i=1}^{m-n} \sum \limits_{j=1, j \neq i}^{m-n}	\hat{s}ec(U_i, U_j) - (m-n)(m-n-1) \| \nabla f \|^2 ,
	\end{align*}
	where $\{ U_1, U_2, \dots, U_{m-n}\}$ is an orthonormal basis of the vertical distribution of $F$. This implies the proof.
\end{proof}

\begin{corollary}
	Let $F: (M^m, g)\to (B^n, g')$ be a Clairaut conformal submersion with $r = e^f$ and dilation $\lambda$ such that the fibers of $F$ are Einstein, with Einstein constant $\lambda_f$. Then 
	\begin{align*}
		Ric(U_i, U_i)& = \lambda_f - (m-n-1) \|\nabla f \|^2,
	\end{align*}
	where $\{U_1, U_2, \dots, U_{m-n}\}$ is an orthonormal basis of the vertical distribution of $F$ and $Ric$ denotes the Ricci curvature of $(M,g)$. 
\end{corollary}
\begin{proof}
	Taking trace of (\ref{eqn3.37}), we get
	\begin{align}\label{eqn3.38}
		\sum \limits_{j=1, j \neq i}^{m-n} sec(U_i, U_j) &= \sum \limits_{j=1, j \neq i}^{m-n}	\hat{s}ec(U_i, U_j) - (m-n-1) \| \nabla f\|^2 ,
	\end{align}
	where $\{ U_1, U_2, \dots, U_{m-n}\}$ is an orthonormal basis of the vertical distribution of $F$, and the sums are taken over all vectors of the orthonormal basis except $U_i$. We know that
	\begin{align}\label{eqn3.39}\small
		\sum \limits_{j=1, j \neq i}^{m-n} sec(U_i, U_j)& = 	\sum \limits_{j=1, j \neq i}^{m-n} \frac{g(R(U_i, U_j)U_j, U_i)}{g(U_i, U_i)g(U_j, U_j) - g(U_i, U_j)^2} \nonumber \\&= Ric(U_i, U_i).
	\end{align}
	Using (\ref{eqn3.39}) in (\ref{eqn3.38}), we get
	\begin{equation}\label{eqn3.40}\small
		Ric(U_i, U_i) = \hat{R}ic(U_i, U_i) - (m-n-1) \|\nabla f \|^2.
	\end{equation}
	Since the fibers of $F$ are Einstein \cite{Besse_1987}, $\hat{R}ic(U_i, U_j) = \lambda_f g(U_i, U_i) = \lambda_f$. Thus (\ref{eqn3.40}) implies the proof.
\end{proof}

\begin{lemma} \label{lemmanormgradf}
	Let $F:(M^m,g) \to (B^n,g')$ be a conformal submersion with constant dilation along fibers, i.e., $\nu \nabla \lambda =0$. Then for every function $f$ such that $\nu \nabla f=0$ we have $\nu \nabla ( g( \nabla f , \nabla f) ) =0$. 
\end{lemma} 
\begin{proof} 
	For every $f$ such that $\nu \nabla f=0$ we can write $f = \phi \circ F$ for some function $\phi$ on $B$. Then for every horizontal $X$ we have
	\begin{equation} \label{eqlemma1}
		g(X, \nabla f) = X(f) = X(\phi \circ F) = (F_* X)(\phi) \circ F = g' ( F_* X , \nabla_B \phi ) \circ F,
	\end{equation}
	where $\nabla_B \phi$ is the gradient of $\phi$ on $(B,g')$. On the other hand, at every point $x \in M$ we have
	\begin{equation} \label{eqlemma2}
		\lambda^2 g(X, \nabla f) = g'(F_* X , F_* \nabla f ) \circ F. 
	\end{equation}
	The two formulas (\ref{eqlemma1}) and (\ref{eqlemma2}) imply that $\nabla f$ projects to $F_* \nabla f = {\tilde \lambda}^2 \nabla_B \phi$, where ${\tilde \lambda}(F(x)) = \lambda(x)$ for all $x \in M$ (which is well defined, because $\nu \nabla \lambda=0$). 
	Hence, the following function on $M$:
	\[
	g' (\nabla_B \phi, \nabla_B \phi) \circ F = \lambda^{-4} ( g'(F_* \nabla f , F_* \nabla f) \circ F) =  \lambda^{-2} g(\nabla f,\nabla f)
	\]
	is constant along fibers, and from $\nu \nabla \lambda =0$ it follows that $\nu \nabla ( g( \nabla f , \nabla f) ) =0$.
\end{proof} 

\begin{corollary}
	If $(M,g)$ is a compact, $m$-dimensional manifold, then every Clairaut conformal submersion from $M$ onto a manifold of dimension $n<m$ has at least two totally geodesic fibers.
\end{corollary}
\begin{proof}
	If $(M,g)$ is compact and $F$ is a Clairaut conformal submersion with $r=e^f$, then there exist points: $x \in M$, where $f$ attains its maximum, and $y \in M$, where $f$ attains its minimum. Then $\nabla f=0$ at $x$ and at $y$. By Theorem \ref{thm3.1} we have $\nu \nabla f=0$, from Lemma \ref{lemmanormgradf} we obtain $\nabla f=0$ along the fiber through $x$ and along the fiber through $y$, and hence by (\ref{eqn3.30}) those fibers are totally geodesic.
\end{proof}

\begin{proposition}
	Let $F: (M^m, g) \to (B^n, g')$ be a Clairaut conformal submersion with $r = e^f$ and dilation $\lambda$. If $(M,g)$ is Einstein, then the metric induced from $g$ on every fiber of $F$ is Einstein.
\end{proposition}

\begin{proof} 
	From \cite[Theorem 2.2.3(1)]{Gudmundsson_1992} and (\ref{eqn3.30}) we obtain for all vertical $V,W$ and an orthonormal basis $\{ U_1, U_2, \dots, U_{m-n}\}$ of the vertical distribution of $F$:
	\begin{align}
		Ric(V,W) &= \hat{R}ic(V, W) + \sum \limits_{i=1}^{m-n} g(T_{U_i} W , T_V U_i ) - \sum \limits_{i=1}^{m-n} g(T_V W , T_{U_i } U_i) \nonumber \\
		&= \hat{R}ic(V, W) + \sum \limits_{i=1}^{m-n} g( {U_i} , W) g(V , U_i ) \| \nabla f \|^2 - \sum \limits_{i=1}^{m-n} g(V , W) g( {U_i } , U_i) \| \nabla f \|^2 \nonumber \\
		&=\hat{R}ic(V, W) - (m-n-1) \| \nabla f \|^2 g( V, W) .
	\end{align}
	Hence, if $(M,g)$ is Einstein, then every fiber of $F$ is Einstein. 
\end{proof}

\section{Harmonicity of Clairaut Conformal Submersions}\label{sec4}
In this section, we discuss the harmonicity of Clairaut conformal submersions.
\begin{definition} 
	(\cite{Sahin_2017}, p. 45) 
	Let $F:(M^m,g) \to (B^n,g')$ be a smooth map between Riemannian manifolds. Then $F$ is harmonic if and only if the tension field $\tau(F)$ of $F$ vanishes at each point $p\in M$.
\end{definition}
\begin{theorem}
	Let $F: (M^m, g)\to (B^n, g')$ be a Clairaut conformal submersion with $r = e^f$ and dilation $\lambda$. Then $F$ is harmonic if and only if
	\begin{equation}\label{harmonicClairaut}
		( m-n )\nabla f	= ( n-2 ) \nabla  \log  { \lambda} . 
	\end{equation}
\end{theorem}
\begin{proof}
	Let $F: (M, g)\to (B, g')$ be a conformal submersion between Riemannian manifolds. Then $\tau(F)= trace (\nabla F_\ast) = \sum_{k=1}^{m}(\nabla F_\ast)(e_k, e_k)= 0$, where $\{e_k\}_{1\leq k \leq m}$ is a local orthonormal basis around a point $p \in M$. Hence
	\begin{equation}\label{eqn3.22}
		\tau(F) = \tau^{kerF_\ast}(F)+ \tau^{(kerF_\ast)^\bot}(F).
	\end{equation}
	Let $\{\tilde{X_j}\}_{1 \leq j \leq n}$,  $\{\lambda X_j\}_{1 \leq j \leq n}$ and   $\{U_i\}_{n+1 \leq i \leq m}$  be orthonormal bases of $TB$, $(kerF_\ast)^\bot$ and $kerF_\ast$, respectively.
	Now, we have \cite{Sahin_2010b}
	\begin{equation}\label{eqn3.24}
		\tau^{kerF_\ast}(F)=-(m-n) F_\ast (H).
	\end{equation}
	Also,
	\begin{equation*}
		\tau^{(kerF_\ast)^\bot}(F)=\sum_{j=1}^{n}(\nabla F_\ast)(\lambda X_j, \lambda X_j).
	\end{equation*}
	Using (\ref{eqn2.11}) in the above equation, we get
	\begin{equation}\label{eqn3.25}
		\tau^{(kerF_\ast)^\bot}(F)= - \frac{\lambda^2}{2}\sum_{j=1}^{n}\left\{2\lambda^2 X_j \left(\frac{1}{\lambda^2}\right)\tilde{X_j} - g(\lambda X_j, \lambda X_j) F_\ast\left(\nabla_{\mathcal{H}}\frac{1}{\lambda^2} \right) \right\}.
	\end{equation}
	Using (\ref{eqn2.12}) in (\ref{eqn3.25}), we get
	\begin{equation*}
		\tau^{(kerF_\ast)^\bot}(F)= - \frac{\lambda^2}{2}\sum_{j=1}^{n}\left\{2\lambda^2 g\left(X_j, \nabla \frac{1}{\lambda^2}\right)\tilde{X_j} - g(\lambda X_j, \lambda X_j) F_\ast\left(\nabla_{\mathcal{H}}\frac{1}{\lambda^2} \right) \right\}.
	\end{equation*}
	Using (\ref{eqn1.1}) in the above equation, we get
	\begin{equation*}
		\tau^{(kerF_\ast)^\bot}(F)= - \frac{\lambda^2}{2}\sum_{j=1}^{n}\left\{2 g'\left(\tilde{X_j}, F_\ast\left(\nabla_{\mathcal{H}}\frac{1}{\lambda^2} \right)\right) \tilde{X_j} - g'(\tilde{X_j}, \tilde{X_j}) F_\ast\left(\nabla_{\mathcal{H}}\frac{1}{\lambda^2} \right) \right\},
	\end{equation*}
	which implies
	\begin{equation}\label{eqn3.26}
		\tau^{(kerF_\ast)^\bot}(F)= (n-2)\frac{\lambda^2}{2}  F_\ast\left(\nabla_{\mathcal{H}}\frac{1}{\lambda^2} \right).
	\end{equation}
	Putting (\ref{eqn3.24}) and (\ref{eqn3.26}) in (\ref{eqn3.22}), we get
	\begin{equation}\label{eqn3.27}
		\tau(F) = (n-2) \frac{\lambda^2}{2} F_\ast\left(\nabla_{\mathcal{H}}\frac{1}{\lambda^2} \right) -(m-n) F_\ast (H).
	\end{equation}
	Since $F$ is a Clairaut conformal submersion, putting $H = -\mathcal{H}\nabla f$ in (\ref{eqn3.27}), we get
	\begin{equation}\label{eqn3.28}
		\tau(F) = (n-2) \frac{\lambda^2}{2} F_\ast\left(\nabla_{\mathcal{H}}\frac{1}{\lambda^2} \right) + (m-n) F_\ast (\mathcal{H}\nabla f).
	\end{equation}
	Hence,
	\[
	\tau(F) = F_\ast {\mathcal{H}} \left(  (n-2) \frac{\lambda^2}{2} \nabla \frac{1}{\lambda^2}  + (m-n) \nabla f \right).
	\]
	As $F_\ast$ is injective on the horizontal distribution, and $\nu \nabla f=0 = \nu \nabla \lambda$, we obtain that vanishing of $\tau(F)$ is equivalent to (\ref{harmonicClairaut}).
\end{proof}

\begin{corollary}
	If $m \neq n=2$ then a Clairaut conformal submersion $F$ is harmonic if and only if it has totally geodesic fibers.
\end{corollary}
\begin{proof}
	For $n=2$ we obtain from (\ref{harmonicClairaut}) that $\nabla f=0$, i.e. the mean curvature of the totally umbilical fibers of $F$ vanishes, and hence those fibers are totally geodesic.
\end{proof}

\begin{corollary}
	Let $F: (M^m, g)\to (B^n, g')$ be a Clairaut conformal submersion with $r = e^f$. If $m \neq 2$, then there exists a function $\phi$ on $M$, constant along the fibers of $F$, such that $F  : (M, \phi^2 g)\to (B, g')$ is a harmonic Clairaut conformal submersion with $r = \phi e^f$.
\end{corollary}
\begin{proof}
	From Lemma \ref{lemmaconformalchangeonM} it follows that $F  : (M, \phi^2 g)\to (B, g')$ is a Clairaut conformal submersion, for which (\ref{harmonicClairaut}) becomes:
	\[
	(m-n)  {\tilde \nabla} (f + \log \phi) =  (n-2)  {\tilde \nabla} \log \left(\frac{ \lambda}{\phi} \right),
	\]
	where ${\tilde \nabla}$ denotes the gradient with respect to $\phi^2 g$. If $m \neq 2$, we can find that the above equation is satisfied if
	\[
	\log \phi =   \frac{n-2}{m-2} \log \lambda - \frac{m-n}{m-2} f + c,
	\]
	for an arbitrary constant $c \in \mathbb{R}$.
\end{proof}

\begin{corollary}
	Let $F: (M^m, g)\to (B^n, g')$ be a Clairaut conformal submersion with $r = e^f$. If $n \neq 2$, then there exists a function $\rho$ on $B$, such that $F  : (M, g)\to (B, \rho^2 g')$ is a harmonic Clairaut conformal submersion with $r = e^f$.
\end{corollary}

\begin{proof}
	From (\ref{eqn1.1}) we obtain for $X,Y \in \mathcal{H}_p$:
	\[
	\rho^2 g'(F_\ast X, F_\ast Y) = \left( (\rho \circ F) \cdot \lambda \right)^2(p) g(X, Y)
	\]
	Hence, $F  : (M, g)\to (B, \rho^2 g')$ is a	Clairaut conformal submersion with $r = e^f$ and dilation
	$(\rho \circ F) \cdot \lambda$, for which equation (\ref{harmonicClairaut}) becomes:
	\[
	( m-n )	\nabla f = ( n-2 ) \nabla (  \log  { \lambda} + \log (\rho \circ F) ) , 
	\]
	which is satisfied for $\rho$ such that
	\[
	\log (\rho \circ F) =  \frac{  m-n }{ n-2  } f - \log \lambda +c,
	\]
	for an arbitrary constant $c \in \mathbb{R}$.
\end{proof}

\section{Clairaut Conformal Submersions with Integrable Horizontal Distributions}\label{sec5}
\begin{example}
	We recall the definition of a doubly warped product $M = M_1 ~{}_\lambda \times_{f_1} M_2$ of two Riemannian manifolds $(M_1, g_1)$ and $(M_2, g_2)$ from \cite{Neill_1983, Ponge_1993}. Let $f_1: M_1 \to \mathbb{R}^+$ and $\lambda: M_2 \to \mathbb{R}^+$ be smooth functions. Let $\pi : M_1 ~{}_\lambda \times_{f_1} M_2 \to M_1$ and $\sigma : M_1 ~{}_\lambda \times_{f_1} M_2 \to M_2$ be the projections. The doubly warped product manifold is the manifold $M = M_1 ~{}_\lambda \times_{f_1} M_2$ furnished with the metric tensor 
	\begin{equation*}
		g= (\lambda \circ \sigma)^2 \pi^\ast(g_1) + (f_1 \circ \pi)^2 \sigma^\ast(g_2),
	\end{equation*}
	that is, if $X$ is tangent to $M_1 ~{}_\lambda \times_{f_1} M_2$ at $(p,q)$, then
	\begin{equation}\label{eqn3.41}
		g(X, X)= \lambda^2 (q) g_1(\pi_\ast X, \pi_\ast X)  + f_1^2(p) g_2(\sigma_\ast X, \sigma_\ast X).
	\end{equation}
	It is easy to prove that the first projection $\pi : M_1 ~{}_\lambda \times_{f_1} M_2 \to M_1$ is a conformal submersion onto $M_1$, whose vertical and horizontal spaces at $(p, q)$ are identified with $T_q M_2$ and $T_p M_1$, respectively, and the second projection  $\sigma : M_1 ~{}_\lambda \times_{f_1} M_2 \to M_2$ is a positive homothetic map onto $M_2$ with scale factor $\frac{1}{f_1^2(p)}$.
	
	Let $\mathcal{L}(M_i)$ for $i=1,2$ be the set of lifts of vector fields from $M_i$ to $M$. We have $U, V \in \mathcal{L}(M_2) \implies [U, V] \in \mathcal{L}(M_2)$, i.e., $\mathcal{H}[U, V] =0$ and $X, Y \in \mathcal{L}(M_1) \implies[X, Y] \in \mathcal{L}(M_1)$, i.e., $\mathcal{\nu}[X, Y] =0$.
	Since $[X, V] =0$, it follows that $[X, V]$ is $\pi$-related to $0$ and $\sigma$-related to $0$. Now, by using Koszul's formula \cite{Neill_1983}, we get
	\begin{equation*}
		\begin{array}{ll}
			g(\mathcal{H} \nabla_U V, Y) &= \frac{1}{2} \{ Ug(V, Y) + V g(Y, U) - Y g(U, V) \\&+ g([U, V], Y) + g([Y, U], V) - g([V, Y], U)\},
		\end{array}
	\end{equation*}
	which implies
	\begin{equation*}
		g(\mathcal{H} \nabla_U V, Y) = -\frac{1}{2}  Y g(U, V).
	\end{equation*}
	By using (\ref{eqn3.41}) in above equation, we get
	\begin{equation}\label{eqn3.43}
		g(\mathcal{H} \nabla_U V, Y) = -\frac{Y}{2}  f_1^2 g_2(\sigma_\ast U, \sigma_\ast V).
	\end{equation}
	Again using (\ref{eqn3.41}) in (\ref{eqn3.43}), we get
	\begin{equation*}
		g(\mathcal{H} \nabla_U V, Y) = -\left(\frac{Yf_1}{f_1}\right) g(U, V).
	\end{equation*}
	Using (\ref{eqn2.3}) in above equation, we get
	\begin{equation*}
		g(T_U V, Y) = -g(Y, \nabla_{\mathcal{H}} \log (f_1)) g(U, V),
	\end{equation*}
	which implies
	\begin{equation}
		T_U V = - g(U, V)\nabla_{\mathcal{H}} \log (f_1),
	\end{equation}
	which means the fibers are totally umbilical. Thus, by Theorem \ref{thm3.1}, $\pi$ is a Clairaut conformal submersion with $r = e^{\log f_1} = f_1$ for any function $\lambda$ constant on the fibers.
\end{example}

Warped product manifold is a local model of the domain of every Clairaut conformal submersion with integrable horizontal distribution.

\begin{proposition}
	Let $F : (M^m,g) \rightarrow (B^n, g')$ be a Clairaut conformal submersion with integrable horizontal distribution, dilation $\lambda$ and $r=e^f$. Then every point of $M$ has a neighbourhood $W$ on which $g$ is a warped product metric.
\end{proposition}

\begin{proof}
	The fibers of $F$ are leaves of the foliation $\mathcal{F}_{\nu}$ on $M$. Since the horizontal distribution of $F$ is integrable, it is tangent to the foliation $\mathcal{F}_{\mathcal{H}}$, which is orthogonal to $\mathcal{F}_{\nu}$. For every $x \in M$ there exists a neighbourhood $W$ of $x$ on which  $\mathcal{F}_{\mathcal{H}}$ defines the submersion $\pi : W \rightarrow {L}$, where ${L}$ is the intersection of $W$ and the leaf of $\mathcal{F}_{\nu}$ passing through $x$; we have $ker\pi_*=(kerF_\ast)^\bot$. We can further assume that $W$ is a tubular neighbourhood of $L$ \cite{Molino}. On $W$, let $U,V  \in \Gamma(kerF_\ast)$ be basic vector fields with respect to the submersion $\pi$ and let $X \in \Gamma(kerF_\ast)^\bot$. Then $\pi_* [X,U] = [\pi_* U , \pi_* X] = [\pi_* U , 0 ] = 0$ and hence $[X,U] \in \Gamma(kerF_\ast)^\bot$. Analogously, we obtain $[X,V] \in \Gamma(kerF_\ast)^\bot$. Since, by Theorem \ref{thm3.1}, the leaves of $\mathcal{F}_{\nu}$ are totally umbilical with $H = - \nabla f$, we have
	\begin{align*}
		-2 g(U,V) X(f) &= 2g(H , X) g(U,V) = g(\nabla_U V , X) + g(\nabla_V U ,X) \\
		& = - g(\nabla_U X , V) - g(\nabla_V X , U) \\
		& = g([ X, U ] , V ) - g(\nabla_X U , V) + g([ X , V ] , U ) - g(\nabla_X V , U) \\
		& = - g(\nabla_X U , V) - g(\nabla_X V , U) \\
		& = - X (g(U,V)).
	\end{align*}
	The above equation has the following solution on $W$, for all $U,V \in kerF_\ast$:
	\begin{equation*}
		g(U,V) = e^{2\psi} g_{ | {L}}(\pi_* U, \pi_* V),
	\end{equation*}
	where $\psi : W \rightarrow \mathbb{R}$ is a function such that $\psi(y) = 0$ for all $y \in L \cap W$ and $\nabla \psi = \nabla f$. Since $\nu \nabla \psi = \nu \nabla f=0$, we can write $\psi = \log (\phi \circ F)$ for some function $\phi$ on $F(W)$, then for all $U,V \in kerF_\ast$ we have
	\begin{equation} \label{gUV}
		g(U,V) =  (\phi \circ F)^2 g_{ | {L}}(\pi_* U, \pi_* V).
	\end{equation}
	Analogously, by (\ref{eqn1.1}) and $\nu \nabla \lambda=0$ (which follows from Theorem \ref{thm3.1}), we can write that for all $X,Y \in (kerF_\ast)^\bot$ we have 
	\begin{equation} \label{gXY}
		g(X,Y) = ({\tilde \lambda} \circ F)^2 g'(F_* X, F_* Y), 
	\end{equation} 
	for a function ${\tilde \lambda}$ on $B$, such that ${\tilde \lambda}(F(x)) = \frac{1}{\lambda(x)}$. From (\ref{gUV}), orthogonality of $kerF_\ast$ and $(kerF_\ast)^\bot$, 
	and (\ref{gXY}), we obtain that the metric $g$ on $W$ is of the form $g = (\phi \circ F)^2 (\pi^* g_{ | {L}} ) + ({\tilde \lambda} \circ F)^2 (F^* g') = (\phi \circ F)^2 (\pi^* g_{ | {L}} ) + F^* ({\tilde \lambda}^2 g')$.
\end{proof}

\begin{example} \label{ex5.2}
	Let $M=\{(x_1,x_2) \in \mathbb{R}^{2} \}$, we define a Riemannian metric $g= e^{2x_2} ( dx_1^2 + dx_2^2 )$ on $M$. Let $B=\{y \in \mathbb{R} \}$ be an Euclidean space with Riemannian metric $g'= dy^2$. Consider a map $F : (M,g) \to (B,g')$ such that
	\begin{equation*}
		F(x_1,x_2)= x_2.
	\end{equation*}
	Then, we get
	\begin{equation*}	
		kerF_\ast = Span\{ U = e_1 \},
	\end{equation*}
	and 
	\begin{equation*}
		(kerF_\ast)^\bot = Span\{ X = e_2\},
	\end{equation*}
	where $\Big\{ e_1 = \frac{\partial}{\partial x_1}, e_2 =\frac{\partial}{\partial x_2}\Big\}$, $\Big\{ e'_1 = \frac{\partial}{\partial y} \Big\}$ are bases on  $T_pM$ and $T_{F (p)}B$ respectively, for all $ p\in M$. We see that $F_\ast (X) = e'_1$ and $\lambda^2 g(X,X)= g'(F_\ast X, F_\ast X)$ for $X \in \Gamma(kerF_\ast)^\bot$ and dilation $\lambda = e^{-x_2}$. Thus $F$ is a conformal submersion.
	
	Now we will find a smooth function $f$ on $M$ satisfying $T_U U =- g(U, U) \nabla f$; $\nu(\nabla \lambda) =0$, for all $U\in \Gamma(kerF_\ast)$. The Christoffel symbols of $g$ are
	\begin{equation*}
		\begin{array}{ll}
			\Gamma_{11}^1 = 0, \Gamma_{11}^2 = -1, \Gamma_{22}^1 = 0, \Gamma_{22}^2 = 1,  \Gamma_{12}^1 = 1= \Gamma_{21}^1, \Gamma_{12}^2 = 0 = \Gamma_{21}^2.
		\end{array}
	\end{equation*}
	By some computations, we obtain
	\begin{equation*}
		\nabla_{e_1} e_1 =  -\frac {\partial}{\partial x_2}, \nabla_{e_2} e_2 =\frac {\partial}{\partial x_2}, \nabla_{e_2} e_1 = \nabla_{e_1} e_2 =\frac {\partial}{\partial x_1}.
	\end{equation*}
	Therefore
	\begin{equation*}
		\nabla_U U = \nabla_{e_1} e_1 = - X.
	\end{equation*}
	Then by (\ref{eqn2.3}), we get
	\begin{equation}\label{eqn3.46}
		T_{U} U= - X.
	\end{equation}
	Since $\lambda= e^{-x_2}$,
	\begin{equation}\label{eqn3.47}
		g\left(U,  \nabla \lambda \right) = 0.
	\end{equation}
	In addition
	\begin{equation}\label{eqn3.48}
		g(U, U) = e^{2x_2}.
	\end{equation}
	For any smooth function $f$ on $M$, the gradient of $f$ with respect to the metric $g$ is given by $\nabla f=\underset{i,j=1}{\overset{2}{\sum}} g^{ij} \frac{\partial f}{\partial x_i} \frac{\partial}{\partial x_j}$. Therefore $\nabla f= e^{-2x_2}\frac{\partial f}{\partial x_1} \frac{\partial}{\partial x_1} +e^{-2x_2} \frac{\partial f}{\partial x_2} \frac{\partial}{\partial x_2}$. Hence, 
	\begin{equation}\label{eqn3.49}
		\nabla f= e^{-2x_2} X,
	\end{equation}
	for the function $f = x_2$. 
	Then, by (\ref{eqn3.46}), (\ref{eqn3.47}), (\ref{eqn3.48}), (\ref{eqn3.49}) and Theorem \ref{thm3.1}, $F$ is a Clairaut conformal submersion.
\end{example}
\vspace{0.3cm}
\noindent \textbf{Acknowledgments} Kiran Meena gratefully acknowledges the research facilities provided by Harish-Chandra Research Institute, Prayagraj, India. In addition, both authors are grateful to the referee for helpful comments that allowed to improve the paper.\\

\noindent \textbf{Authors Contributions} Conceptualization, methodology, investigation, writing - original draft, validation, review, editing and reading have been performed by both the authors and agreed to the paper.\\

\noindent \textbf{Funding} Kiran Meena is financial supported by the Department of Atomic Energy, Government of India [Offer Letter No.: HRI/133/1436 Dated 29 November 2022].\\

\noindent \textbf{Data Availability} This manuscript has no associated data.\\

\noindent \textbf{Code Availability} Not Applicable.

\section*{Declarations}
\vspace{0.3cm}
\noindent \textbf{Conflict of interest} The authors have no Conflict of interest and no financial interests for this article.\\

\noindent \textbf{Ethical Approval} The submitted work is original and not submitted to more than one journal for simultaneous consideration.\\

\noindent \textbf{Consent to Participate} Not applicable.\\

\noindent \textbf{Consent for Publication} Not applicable.\\

\noindent \textbf{Open Access} This article is licensed under a Creative Commons Attribution 4.0 International License, which permits use, sharing, adaptation, distribution and reproduction in any medium or format, as long as you give appropriate credit to the original author(s) and the source, provide a link to the Creative Commons licence, and indicate if changes were made. The images or other third party material in this article are included in the article's Creative Commons licence, unless indicated otherwise in a credit line to the material. If material is not included in the article’s Creative Commons licence and your intended use is not permitted by statutory regulation or exceeds the permitted use, you will need to obtain permission directly from the copyright holder. To view a copy of this licence, visit \textcolor{blue}{http://creativecommons.org/licenses/by/4.0/}.\\

\vspace{1cm}
\noindent \textbf{Publisher's Note} Springer Nature remains neutral with regard to jurisdictional claims in published maps and institutional affiliations.
\end{document}